\documentclass[journal]{IEEEtran}

\usepackage{mathrsfs,amsmath,amssymb,amsthm,graphicx,xcolor,dsfont}
\usepackage{multirow,soul,framed}
\usepackage[comma,numbers,square,sort&compress]{natbib}

\colorlet{shadecolor}{cyan!20}

\newtheorem{assumption}{Assumption}

\newtheorem{theorem}{Theorem}
\newtheorem{remark}{Remark}

\newtheorem{corollary}{Corollary}

\newcommand{\R}{\mathbb{R}}

\usepackage[prependcaption,colorinlistoftodos]{todonotes}

\begin{document}
\title{Data-driven stabilization of nonlinear polynomial systems with noisy data}

\author{Meichen Guo, Claudio De Persis, and Pietro Tesi% <-this % stops a space
%\thanks{This work was not supported by any organization}% <-this % stops a space
\thanks{Meichen Guo and Claudio De Persis are with ENTEG, Faculty of Science and Engineering, University of Groningen, 9747 AG Groningen, The Netherlands. Email: {\tt\small meichen.guo@rug.nl, c.de.persis@rug.nl} Pietro Tesi is with DINFO, University of Florence, 50139 Florence, Italy. E-mail: {\tt\small pietro.tesi@unifi.it}}
}

\maketitle
\date{}

\begin{abstract}
In a recent paper we have shown how to learn controllers for unknown linear systems using finite-sized noisy data by solving linear matrix inequalities. In this note we extend this approach to deal with unknown nonlinear polynomial systems by formulating stability certificates in the form of data-dependent sum of squares programs, whose solution directly provides a stabilizing controller and a Lyapunov function. We then derive variations of this result that lead to more advantageous controller designs. The results also reveal connections to the problem of designing a controller starting from a least-square estimate of the polynomial system.
%a method to robustly design controllers starting from a least-square estimate of the unknown polynomial system.
\end{abstract}

\section{Introduction}\label{section:Introduction}

The idea of designing control laws from measured data alone for a system with unknown dynamics is
of high relevance, both theoretically and practically. Following a rough classification, there are two approaches to concretize this idea. The first one is based on identifying a model of the system and then designing a  controller that deals with the uncertainty  resulting from the identification process. The second one devises methods to directly synthesize controllers  from data without explicitly undertaking a model identification, hence the name of direct data-driven control.

Both approaches are viable and appealing but
%have pros and cons -- see  \cite{sznaier.position} for a modern excursus on these topics --
the focus of this paper is on direct data-driven control, which will lead us to provide neat  compact
conditions for the synthesis of nonlinear stabilizing controllers and corresponding Lyapunov functions. Various results on direct data-driven  control have been proposed and much attention has been devoted to the specific study of nonlinear systems, which remains an open and challenging problem. Among the contributions on direct data-driven control of nonlinear systems, we recall works such as the virtual reference feedback tuning (VRFT) \cite{campi2006VRFT} and iterative feedback tuning \cite{Hjalmarsson1998IFT},
the on-line data-driven control of \cite{tanaskovic2017data},
the so-called intelligent PID \cite{Fliess2009,Tabuada2017CDC}, and reinforcement learning \cite{Lee2005nonlinearQlearning,Lopez2019RL}.

{\it Related work.} Recent works inspired by Willems et al.'s lemma \cite{Willems2005} have revealed a new approach to accomplish direct data-driven control. Specifically, for linear systems, by Willems et al.'s lemma, any input-output trajectory of a linear system can be expressed as a linear combination of finite number of previously measured input-output data. Using this lemma to parameterize a closed-loop linear system controlled by a state-feedback controller, the feedback gain can be directly found via data-based linear matrix inequalities (LMIs) for various control problems \cite{cdpTAC2020}. A nonlinear extension of the linear stabilization result was also considered in \cite{cdpTAC2020}, where the nonlinear remainder after linearization is treated as disturbances and robust data-driven control method is applied to obtain local stabilization result. Recently, \cite{Bisoffi2020bilinear} investigated the stabilization of bilinear systems with a characterization of the basin of attraction. Polynomial systems lend themselves to a similar analysis as the one in \cite{cdpTAC2020} once one expresses the system's polynomial vector fields in the span of a suitable basis of polynomial functions and pursues a simultaneous design of controllers and Lyapunov functions based on the second Lyapunov theorem \cite{guoCDC2020}. Independently, and taking a different route,
\cite{dai2020semialgebraic} used
the dual stability theory
and Farkas' lemma for direct data-driven design of rational state-feedback controllers for nonlinear systems. A basis of polynomials to express the nonlinearity was considered earlier in \cite{novara16data-driven}
and used to perform a  right-inversion of the system dynamics based on which a reference tracking controller is derived. Model-inversion errors are then compensated by an outer-loop linear controller designed according to the VRFT technique.

A reasonable concern in data-driven control is that the measured data can be affected by unknown noise. For linear systems, \cite{cdpTAC2020} posed a signal-to-noise ratio assumption on the measurement noise and presented sufficient conditions for data-driven stabilizers. This method was further investigated and analyzed in \cite{cdpLQRnoise} for data-driven design of linear quadratic regulators. The paper \cite{Berberich2019robust} has remarked that the condition on the noise introduced in \cite{cdpTAC2020} can be interpreted as an instance of a quadratic matrix inequality and  analyzed via a full-block S-procedure. Another work \cite{vanwaarde2020noisy} presented a matrix-valued S-lemma based on the classical S-lemma and apply it to data-driven control of linear systems. Other notable works dealing with direct or indirect data-driven robust control include, but are not limited to, \cite{Dai2019cdc,dai2020semialgebraic}.%,Recht2019LQR_CE,Wabersich2020}.

In this paper, we investigate global data-driven stabilization of continuous-time nonlinear polynomial systems using noisy data. The focus is on polynomial systems because, first, polynomials systems are widely used to model processes in engineering applications such as fluid dynamics \cite{Chernyshenko2014Fluid} and robotics \cite{Majumdar2013ICRA}.
In fact, polynomial control systems, and the supporting technical developments in sum of squares (SOS) optimization, have attracted considerable attention over the last twenty years \cite{ParriloThesis,henrion2005positive,chesi2009tac,Chesi10LMI,Hancock2013AUT,Valmorbida2017aut,Ahmadi2019SIAAG}.
Second, as shown in the model-based control \cite{Prajna2004SOS}, nonlinear polynomial systems can be written into a linear-like form and controlled using the Lyapunov method. This inspires us to adopt the framework of \cite{cdpTAC2020} to establish a direct data-driven control synthesis for the nonlinear polynomial systems. Our previous work \cite{guoCDC2020} has presented some preliminary results on data-driven stabilization of polynomial systems with state-independent input vector field and noise-free data. In this paper, we enhance the results to handle state-dependent input vector field and noisy data. Besides being of interest in its own right, the study on polynomial systems will help us better understand and gain more insights in direct data-driven control of nonlinear systems.

{\it Contribution.} The main contribution of this paper is developing state feedback stabilizers for nonlinear polynomial systems using noisy input-state data alone. It is assumed that the noise is unknown and it satisfies a quadratic constraint. Using a variation of the data-based closed-loop representation in our work \cite{guoCDC2020}, we improve the previous results so that high-order polynomial systems can now be handled more effectively. Following the framework of \cite{cdpTAC2020}, the first result of this paper shows that the data-driven design with measurement noise for linear systems can be extended to nonlinear polynomial systems using the adjusted closed-loop representation and matrix-valued Young's inequality. Then, by adopting other forms of the data-based closed-loop representation, more advantageous stabilization results are attained which require less assumptions or have improved computational efficiency. A connection is also established with the problem of designing controllers based on least-square estimates of the system's dynamics. As the data-driven stabilizing conditions involve semi-positive-definiteness of matrix polynomials, we use SOS relaxations to make the computation tractable. We also show that the results can be interpreted as SOS relaxations of pointwise necessary and sufficient conditions for global robust stabilizability of
polynomial systems.

{\it Organization.} The structure of the paper is specified as follows. In Section \ref{section:preliminaries}, we formulate the data-driven control problem. Section \ref{section:mainresults} contains the main results where stabilizers are designed using data corrupted by measurement noise. Theorem \ref{noisy-stabilization-I} can be regarded as an extension of \cite[Theorem 5]{cdpTAC2020} from linear systems to nonlinear polynomial systems. Next, by slightly modifying the closed-loop system representation, a result in the same spirit of Theorem \ref{noisy-stabilization-I} is proposed that requires less assumptions on the controlled dynamics. A more computationally efficient corollary then follows, where a decision variable is redefined such that its size is independent of the size of the data. In Section \ref{section:examples}, the simulation results on the Van der Pol oscillators are presented to show the effectiveness of the proposed design methods. Finally, we summarize the paper and draw the conclusions in Section \ref{section:conclusion}.

{\it Notation.} The following notations are adopted throughout the paper:
\begin{itemize}
  \item[--] $A\succeq 0$: matrix $A$ is positive semi-definite;
  \item[--] $A\succ 0$: matrix $A$ is positive definite;
  \item[--] $A\succeq B$: matrix $A-B$ is positive semi-definite;
  %\item[--] $A\otimes B$: Kronecker product of matrices $A$ and $B$;
  %\item[--] $\mathbb{N}$: the set of natural numbers including $0$;
  \item[--] $\mathbb{N}_{>0}$: the set of natural numbers excluding $0$;
  \item[--] $\R$: the set of real numbers;
  \item[--] $\R_{>0}$: the set of positive real numbers;
  %\item[--] $\mathbb{S}^{n}$: the set of $n \times n$ symmetric matrices;
  \item[--] $\mathcal{P}$: the set of polynomials;
  \item[--] $\Sigma$: the set of SOS polynomials;
  \item[--] $\mathcal{P}^{r\times s}$: the set of $s\times r$ matrix polynomials;
  \item[--] $\Sigma^{r}$: the set of $r\times r$ SOS matrix polynomials.
  %\item[--] $\partial\mathbb{X}$: the boundary of a compact set $\mathbb{X}$.
\end{itemize}

\section{Problem Formulation}
\label{section:preliminaries}

Consider the polynomial system
\begin{align}\label{polysyst:original}
\dot{x} =  f(x)+g(x)u
\end{align}
where $x\in\R^{n}$ is the state, $u\in\R^{m}$ is the control input, $f(x)$ and $g(x)$ are polynomial vector fields of sizes $n \times 1$  and $n\times m$, respectively. The specific expressions of $f(x)$ and $g(x)$ are unknown. The polynomial system \eqref{polysyst:original} can be written into the linear-like form
\begin{align}\label{polysyst:linearlike}
\dot{x} = AZ(x)+BW(x)u
\end{align}
where $A\in\mathbb{R}^{n\times N}$ and $B\in\mathbb{R}^{n\times q}$ are unknown constant matrices.
The $N\times 1$ vector $Z(x)$ is a collection of distinct monomials in $x$ that may appear in $f(x)$, and the $q\times m$ matrix $W(x)$ contains monomials that may appear in $g(x)$.
\medskip

%\cmargin{I would omit this. You discuss it after \eqref{Z-hatZ}}
%\cdp{We need the following assumption on the polynomial system \eqref{polysyst:original}.
%\begin{assumption}\label{assumption:equilibrium}
%The autonomous system $\dot{x}=f(x)$ has an equilibrium at the origin, \emph{i.e.}, $Z(x)=0$ if $x=0$.
%\end{assumption}
%}

As the dynamics of \eqref{polysyst:original} is unknown, to write the linear-like form \eqref{polysyst:linearlike}, we need the vector $Z(x)$ to contain all distinct monomials that may appear in $f(x)$. For nonlinear systems having high order, $Z(x)$ can have notably large size and high order, which may cause issues in the controller design.
As pointed out in \cite{Prajna2004SOS} for the model-based control of polynomial systems, the choice of $Z(x)$ affects the success of the design method. A big vector $Z(x)$ containing high order monomials can cause computational issues in the SOS program and fail to give a solution. To overcome this issue, we use another vector $\hat{Z}(x)$ having smaller size and lower degree than $Z(x)$ for the controller design and stability analysis. Specifically, the $p\times 1$ vector $\hat{Z}(x)$ satisfies that
\begin{align}
Z(x) = H(x)\hat Z(x)\label{Z-hatZ}
\end{align}
%where
with matrix polynomial $H(x)\in\mathcal{P}^{N\times p}$. Our controller design and stability analysis will be based on the Lyapunov function $V(x)=\hat{Z}(x)^{\top}P^{-1}\hat{Z}(x)$, with $P$ a positive definite matrix. Since in this paper we restrict ourselves to consider a global stabilization problem,
we have that $\hat{Z}(x)=0$ if and only if $x=0$, and $\hat{Z}(x)$ is radially unbounded. A straightforward choice satisfying these conditions is a $\hat{Z}(x)$ whose first $n$ components coincide with $x$. This choice will be used in the subsequent sections of this paper. Note that $\hat{Z}(x)=0$ if and only if $x=0$ implies that $Z(x)=0$ if $x=0$, \emph{i.e.}, the autonomous system $\dot{x}=f(x)$ must have an equilibrium at the origin. The latter results in no loss of generality because, in the case $f(x)$ has no equilibrium at the origin,
%and the equilibrium to which the system must be stabilized is known,
one reduces the analysis to the case $f(0)=0$ by a change of state variables.
%{\color{violet}
On the other hand, even though the matrix $W(x)$ having larger size also increases computational burden in the SOS program, it does not directly affect the design of the Lyapunov function or the sizes of the decision variables. Hence, the choice of $W(x)$ has much less impact on the success of the synthesis and will not be a focus of this work.
%}

To design a data-driven controller, we first run an off-line experiment to collect the input-state data. The experiment is conducted over the time interval $[t_{0},t_{0}+(T-1)\tau]$ where $T\in\mathbb{N}_{>0}$ is the number of sampled data and $\tau\in\R_{>0}$ is the sampling time. The Hankel matrices of the sampled input-state %output
data are defined as
\begin{align*}%\label{}
  U_{0}&%=U_{0,1,T}
  :=\begin{bmatrix}
              u(t_{0})\! & \!u(t_{0}\!+\!\tau)\! & \!\cdots\! &\! u(t_{0}+(T\!-\!1)\tau)
            \end{bmatrix},\\
  X_{0}&
  :=\begin{bmatrix}
              x(t_{0})\! & \!x(t_{0}\!+\!\tau)\! & \!\cdots\! &\! x(t_{0}+(T\!-\!1)\tau)
            \end{bmatrix},\\
  X_{1}&
  :=\begin{bmatrix}
              \dot{x}(t_{0})\! & \!\dot{x}(t_{0}\!+\!\tau)\! & \!\cdots\! &\! \dot{x}(t_{0}+(T\!-\!1)\tau)
            \end{bmatrix}.
\end{align*}
Then, we can compute the data
\begin{align*}
  Z_{0}
  &:=\!\big[ Z(x(t_{0}))~~ Z(x(t_{0}\!+\!\tau)) ~~ \cdots ~~ Z(x(t_{0}\!+\!(T\!-\!1)\tau))\big] ,\\
 \overline{U}_{0}
  &:=\!\big[ W(x(t_{0}))u(t_{0}) ~~ W(x(t_{0}\!+\!\tau))u(t_{0}\!+\!\tau)\\
&\quad\quad\quad\quad~~ \cdots ~~ W(x(t_{0}+(T\!-\!1)\tau))u(t_{0}+(T\!-\!1)\tau)\big].
\end{align*}

We consider experiments affected by measurement noises $D_{0}$ defined as
\begin{align*}
D_{0}%=D_{0,T}
:=\begin{bmatrix}
              d(t_{0})\! & \!d(t_{0}\!+\!\tau)\! & \!\cdots\! &\! d(t_{0}+(T\!-\!1)\tau)
            \end{bmatrix}
\end{align*}
where $d(t)\in\mathbb{R}^{n}$. The specific nature of the noise is discussed in the next section. Affected by the noise, the measured derivatives $X_{1}$ satisfies
\begin{align}
X_{1} = AZ_{0} + B\overline{U}_{0} + D_{0}. \label{data-equation}
\end{align}

The data-driven stabilization problem of \eqref{polysyst:original} is to design a state-dependent control gain $F(x)$ using the noisy experimental input-state data alone, such that under the state-feedback controller
\begin{align}
u =  F(x)\hat{Z}(x)
\end{align}
the closed-loop system is globally asymptotically stable at the origin.

\section{Data-driven stabilization with noisy data}
\label{section:mainresults}

In this section, we present conditions for the stabilization of unknown polynomial systems with data affected by noise.
%Consider experiments affected by measurement noises $D_{0}$ defined as
%\begin{align*}
%D_{0}%=D_{0,T}
%:=\begin{bmatrix}
%              d(t_{0})\! & \!d(t_{0}\!+\!\tau)\! & \!\cdots\! &\! d(t_{0}+(T\!-\!1)\tau)
%            \end{bmatrix}
%\end{align*}
%where $d(t)\in\mathbb{R}^{n}$. Affected by the noise, the measured derivatives $X_{1}$ satisfies
%\begin{align}
%X_{1} = AZ_{0} + B\overline{U}_{0} + D_{0}. \label{data-equation}
%\end{align}
%
It is unrealistic to expect to stabilize the system using data that are affected by arbitrary noise. Hence, we introduce the following assumption on $D_{0}$ :
\begin{assumption}\label{assumption:boundaryD}
The matrix $D_{0}$ satisfies $D_{0}D_{0}^{\top}\preceq R_{D}R_{D}^{\top}$ for some known $R_{D}\in\mathbb{R}^{n\times T}$.
\end{assumption}

Assumption \ref{assumption:boundaryD} is in a similar form as \cite[Assumption 2]{cdpTAC2020} where a signal-to-noise ratio assumption is posed on the noise $D_{0}$ to solve the data-driven control problems. In fact, if we set $R_{D}=\gamma X_{1}$ for some constant $\gamma>0$, Assumption \ref{assumption:boundaryD} becomes the same as \cite[Assumption 2]{cdpTAC2020}. Interestingly, the noisy data result in \cite[Theorem 5]{cdpTAC2020} can be extended to nonlinear polynomial systems, which will be presented in Subsection \ref{subsection:noisy1}. Then, in Subsection \ref{subsection:noisy2} and \ref{subsection:noisy3}, by slightly changing the data-based closed-loop representation and the decision variable, we can derive the result under more advantageous conditions.

\subsection{Data-driven stabilization with noisy data}
\label{subsection:noisy1}

In the result below, inspired by the results of \cite{cdpTAC2020}, we will introduce a controller of the form (a) $u=F(x)\hat Z(x)$ where $F(x)=U_{0}Y(x)P^{-1}$ and  $Y(x), P \succ 0$ are matrices
that satisfy the condition (b) $Z_{0}Y(x) = H(x)P$. This choice allows us to express the closed-loop dynamics as
\begin{equation}\label{noisy-closed-loop}
\begin{array}{rcl}
\dot{x} & = & AZ(x)+BW(x)u \\
& %\stackrel{\eqref{control:thm2}}{=}
\stackrel{{\rm (a)}}{=}
&\big( AH(x)P P^{-1}+BW(x)U_{0}Y(x)P^{-1} \big)\hat{Z}(x)\\
& %\stackrel{\eqref{condition:thm2-1}}{=}
\stackrel{{\rm (b)}}{=}
&\big( AZ_{0} Y(x)P^{-1}+BW(x)U_{0} Y(x)P^{-1} \big)\hat{Z}(x)\\
& = &(AZ_{0}+ BW(x)U_{0})Y(x)P^{-1} \hat{Z}(x)\\
& \stackrel{\eqref{data-equation}}{=} &(X_{1} -D_0-B \overline{U}_{0} + BW(x)U_{0}) Y(x)P^{-1} \hat{Z}(x)
\\
& = & (X_{1} -E \hat{U}_{0}(x)) Y(x)P^{-1} \hat{Z}(x)
\end{array}
\end{equation}
having introduced the symbols
\begin{align*}
E := [B ~~D_{0}], \quad
\widehat{U}_{0}(x):= \begin{bmatrix}\,\overline{U}_{0}-W(x)U_{0}\\ I_{T}\end{bmatrix}
\end{align*}
to obtain the last identity.

To design the stabilizing gain $F(x)$, we need another assumption on the unknown matrix $B$.
\begin{assumption}\label{assumption:boundaryB}
The matrix $B$ satisfies $BB^{\top}\preceq R_{B}R_{B}^{\top}$ for some known $R_{B}\in\mathbb{R}^{n\times q}$.
\end{assumption}

Under Assumptions \ref{assumption:boundaryD} and \ref{assumption:boundaryB},
\[
EE^{\top} = BB^{\top} + D_{0}D_{0}^{\top}\preceq R_{B}R_{B}^{\top}+R_{D}R_{D}^{\top}.
\]
Define $R_{E}:=[R_{B} ~~R_{D}]$. Then, it holds that
\begin{equation}\label{cond:boundaryE}
EE^{\top} \preceq  R_{E}R_{E}^{\top}.
\end{equation}

Obtaining stabilization results requires to impose positive conditions on matrix polynomials, which is often computationally intractable. To resolve this issue, we use SOS relaxations in the subsequent results.

\begin{theorem}\label{noisy-stabilization-I}
For the polynomial system (\ref{polysyst:original}), under Assumptions \ref{assumption:boundaryD} and \ref{assumption:boundaryB}, if there exist a positive definite matrix
$P\in\mathbb{R}^{p\times p}$, matrix polynomial $Y(x)\in\mathcal{P}^{T\times p}$, $\epsilon_{1}(x)\in\Sigma$, and $\epsilon_{2}>0$
such that
\begin{align}
&Z_{0}Y(x) = H(x)P,\label{condition:I-1}\\
%&P-\mu I_{p} \succeq 0,\label{condition:I-2}\\
& \begin{bmatrix}
\Upsilon_{E}(x)- \epsilon_{1}(x)I_{p} & Y(x)^{\top}\widehat{U}_{0}(x)^{\top}\\
\widehat{U}_{0}(x)Y(x) & \epsilon_{2} I_{q+T}
\end{bmatrix} \in\Sigma^{p+q+T},\label{condition:I-3}
\end{align}
where %$\epsilon_{1}(x)\in\Sigma$,
%$\epsilon_{2}(x)>0$ for all $x\in\mathbb{R}^{n}$,
%$\epsilon_{2}>0$,
%and
\begin{align*}
\Upsilon_{E}(x) &:= -\frac{\partial\hat{Z}(x)}{\partial x}X_{1}Y(x) - Y(x)^{\top}X_{1}^{\top}\frac{\partial\hat{Z}(x)}{\partial x}^{\top}\\
&\quad~ -\epsilon_{2}
\frac{\partial\hat{Z}(x)}{\partial x} R_{E} R_{E} ^{\top} \frac{\partial\hat{Z}(x)}{\partial x}^{\top},
\end{align*}
then the state-feedback controller
\begin{align}\label{control:thm1}
u = U_{0}Y(x)P^{-1}\hat{Z}(x)
\end{align}
%globally
stabilizes the continuous-time polynomial system \eqref{polysyst:original}. Moreover, if $\epsilon_{1}(x)>0$ for all $x\ne 0$, the closed-loop system is globally asymptotically stable.
\end{theorem}

\proof

Design the Lyapunov function as $V(x)=\hat{Z}(x)^{\top}P^{-1}\hat{Z}(x)$ whose derivative along the dynamics \eqref{noisy-closed-loop} is
\begin{align}\label{lyap.deriv.id}
\dot{V}(x) = -\hat{Z}&(x)^{\top} P^{-1} \Big[ -\frac{\partial\hat{Z}(x)}{\partial x}\big( X_{1}\!-\!E\widehat{U}_{0}(x) \big)Y(x)\nonumber\\
&- Y(x)^{\top}\big( X_{1}\!-\!E\widehat{U}_{0}(x) \big)^{\top} \frac{\partial\hat{Z}(x)}{\partial x}^{\top} \Big] P^{-1}\hat{Z}(x).
\end{align}
The SOS condition (\ref{condition:I-3}) guarantees the non-negativity of the matrix, which, by Schur complement, is equivalent to
\begin{align*}
\Upsilon_{E}(x)- \epsilon_{2}^{-1}
Y(x)^{\top}\widehat{U}_{0}(x)^{\top} \widehat{U}_{0}(x)Y(x) \succeq \epsilon_{1}(x)I_{p}
\end{align*}
for all $x\in\mathbb{R}^{n}$. Under Assumptions \ref{assumption:boundaryD} and \ref{assumption:boundaryB}, %{\color{olive} The general form of Young's inequality is $L(x)M(x)^{\top}+M(x)L(x)^{\top} \preceq \epsilon_{2}(x)L(x)L(x)^{\top} + \epsilon_{2}(x)^{-1}M(x)M(x)^{\top}$}
the left-hand side of the inequality satisfies
\begin{align*}
&\quad \Upsilon_{E}(x)- \epsilon_{2}^{-1}
Y(x)^{\top}\widehat{U}_{0}(x)^{\top} \widehat{U}_{0}(x)Y(x)\\
%&=-\frac{\partial\hat{Z}(x)}{\partial x}X_{1}Y(x) - Y(x)^{\top}X_{1}^{\top}\frac{\partial\hat{Z}(x)}{\partial x}^{\top}\\
%&\quad -\epsilon_{2}(x) \frac{\partial\hat{Z}(x)}{\partial x} R_{E} R_{E} ^{\top} \frac{\partial\hat{Z}(x)}{\partial x}^{\top}\\
%&\quad - \epsilon_{2}(x)^{-1}Y(x)^{\top}\widehat{U}_{0}(x)^{\top} \widehat{U}_{0}(x)Y(x)\\
&\preceq -\frac{\partial\hat{Z}(x)}{\partial x}X_{1}Y(x) - Y(x)^{\top}X_{1}^{\top}\frac{\partial\hat{Z}(x)}{\partial x}^{\top}\\
&\quad~-\epsilon_{2}
\frac{\partial\hat{Z}(x)}{\partial x} EE ^{\top} \frac{\partial\hat{Z}(x)}{\partial x}^{\top}\\
&\quad~-\epsilon_{2}^{-1}
Y(x)^{\top}\widehat{U}_{0}(x)^{\top} \widehat{U}_{0}(x)Y(x)\\
&\preceq -\frac{\partial\hat{Z}(x)}{\partial x}X_{1}Y(x) - Y(x)^{\top}X_{1}^{\top}\frac{\partial\hat{Z}(x)}{\partial x}^{\top}\\
&\quad~+ \frac{\partial\hat{Z}(x)}{\partial x} E \widehat{U}_{0}(x)Y(x)
+Y(x)^{\top}\widehat{U}_{0}(x)^{\top}E ^{\top} \frac{\partial\hat{Z}(x)}{\partial x}^{\top}\\
&=-\frac{\partial\hat{Z}(x)}{\partial x}\big( X_{1}-E\widehat{U}_{0}(x) \big)Y(x)\\
&\quad~- Y(x)^{\top}\big( X_{1}-E\widehat{U}_{0}(x) \big)^{\top} \frac{\partial\hat{Z}(x)}{\partial x}^{\top},
\end{align*}
where the last inequality follows from the matrix-valued Young's inequality
$\epsilon_{2} L(x)L(x)^{\top}+\epsilon_{2}^{-1}
M(x)M(x)^{\top}\succeq -L(x)M(x)^{\top}-M(x)L(x)^{\top} $.
Thus, it holds that
\begin{align*}
&-\!\frac{\partial\hat{Z}(x)}{\partial x}\big( X_{1}\!-\!E\widehat{U}_{0}(x) \big)Y(x)\\
&\quad~ - Y(x)^{\top}\big( X_{1}\!-\!E\widehat{U}_{0}(x) \big)^{\top} \frac{\partial\hat{Z}(x)}{\partial x}^{\top} \succeq \epsilon_{1}(x)I_{p}
\end{align*}
for all $x\in\mathbb{R}^{n}$. Then, the derivative of $V(x)$ satisfies that
\begin{align}\label{cond:dotV-I}
\dot{V}(x) \le -\epsilon_{1}(x) \hat{Z}(x)^{\top} P^{-1} \cdot P^{-1}\hat{Z}(x)
\end{align}
for all $x\in\mathbb{R}^{n}$. Therefore, the closed-loop system is
%globally
stable. Moreover, if $\epsilon_{1}(x)>0$ for all $x\ne 0$, then $\dot{V}(x)<0$ for all $x\ne 0$, and thus the closed-loop system is globally asymptotically stable at the origin since $V(x)$ is radially unbounded (recall that $\hat Z(x)$ contains $x$ as a sub-vector). \qed
\medskip

To find $P$ and $Y(x)$ in Theorem \ref{noisy-stabilization-I}, one can use MATLAB based tools such as the SOSTOOLS \cite{papachristodoulou2018sostools}. The decision variables of the SOS program include the constant matrix $P$ and the coefficients of $Y(x)$. The coefficients of $\epsilon_{1}(x)$ and
$\epsilon_{2}$ can also be set as decision variables.

\begin{remark}
\emph{(Input vector field independent of $x$)}\label{rem:inp-state-ind}
In the case where the input vector field is independent of $x$, Theorem \ref{noisy-stabilization-I} can achieve data-driven stabilization without Assumption \ref{assumption:boundaryB}. In this case, the polynomial system can be written into the linear-like form
\begin{align*}
\dot{x} =  AZ(x)+Bu
\end{align*}
where $B\in\mathbb{R}^{n\times m}$, and the corresponding closed-loop dynamics can be expressed using data as
\begin{align*}
\dot{x} = (X_{1}-D_0)Y(x)P^{-1}\hat{Z}(x).
\end{align*}
Then, a stabilization result in the same fashion as Theorem \ref{noisy-stabilization-I} can be derived by changing \eqref{condition:I-3} into
\begin{align*}
\begin{bmatrix}
\Upsilon(x)- \epsilon_{1}(x)I_{p} & Y(x)^{\top}\\
Y(x) & \epsilon_{2}(x)I_{T}
\end{bmatrix} \in\Sigma^{p+T}
\end{align*}
where
\begin{align*}
\Upsilon(x) &:= -\frac{\partial\hat{Z}(x)}{\partial x}X_{1}Y(x) - Y(x)^{\top}X_{1}^{\top}\frac{\partial\hat{Z}(x)}{\partial x}^{\top}\\
&\quad~ -\epsilon_{2}(x) \frac{\partial\hat{Z}(x)}{\partial x} R_{D} R_{D} ^{\top} \frac{\partial\hat{Z}(x)}{\partial x}^{\top}.
\end{align*}
\end{remark}

\begin{remark}\label{rem:rank_Z_0}
\emph{(Rank condition on $Z_{0}$)} The feasibility of condition \eqref{condition:I-1} implies a rank condition of $Z_{0}$. As $P\succ 0$, $\text{rank}\big(H(x)P\big)=\text{rank} \big( H(x)\big)$. Based on the assumptions on $Z(x)$ and $\hat{Z}(x)$, the rank of $H(x)\in \mathcal{P}^{N\times p}$ should satisfy $1\le \text{rank}\big( H(x)\big)\le p$. This means that $Z_{0}Y(x)$ should also have the same rank between $1$ and $p$. Hence, the rank condition needed for $Z_{0}$ is $\text{rank}\big( Z_0 \big) \ge \text{rank}\big( H(x)\big)$. Depending on $H(x)$, in many practical cases $Z_{0}$  must have  full row rank in order to fulfil condition \eqref{condition:I-1}. These conditions on $Z_{0}$ imply that the number $T$ of sampling data must be sufficiently large, since $Z_{0}\in \mathbb{R}^{N\times T}$.
\end{remark}

\subsection{Data-driven stabilization without bounds on the input matrix $B$}
\label{subsection:noisy2}

By defining $S:=[B~~A]$ and $W_{0}(x):=[U_{0}^{\top}W(x)^{\top} ~~ Z_{0}^{\top}]^{\top}$, the chain of equalities in  \eqref{noisy-closed-loop} has highlighted the identity %$X_{1} -E \hat{U}_{0}(x)= SW_{0}$
\begin{align}\label{a-key-identity}
X_{1} -E \hat{U}_{0}(x)= SW_{0}(x).
\end{align}
%The data-driven stabilization result in Proposition \ref{noisy-stabilization-I} requires Assumption \ref{assumption:boundaryB} on the unknown matrix $B$. In what follows, we slightly change the data-based representation by putting the unknown dynamics in one matrix $S:=[B~~A]$. In doing so, we present a similar stability result which contain a few advantages compared with Proposition \ref{noisy-stabilization-I}.
%
%
%In view of the identity \eqref{a-key-identity},
We can expect then that the analysis carried out in the previous subsection for  the representation
$\dot x= (X_{1} -E \hat{U}_{0}(x)) Y(x)P^{-1} \hat{Z}(x)$ can be repeated for the representation
\begin{align}\label{closdlp:S}
\dot x = SW_{0}(x)Y(x)P^{-1}\hat{Z}(x)
\end{align}
and the same stability result can be established.
%The data-driven stabilization result in Proposition \ref{noisy-stabilization-I} requires Assumption \ref{assumption:boundaryB} on the unknown matrix $B$. In what follows, we slightly change the data-based representation by putting the unknown dynamics in one matrix $S:=[B~~A]$.
In doing so, we obtain a few advantages compared to Theorem \ref{noisy-stabilization-I}, including the relaxation of Assumption \ref{assumption:boundaryB} on the unknown matrix $B$.

As a first remark, standard calculations show that the derivative of the Lyapunov function $V(x)=\hat{Z}(x)^{\top}P^{-1}\hat{Z}(x)$ along the closed-loop system $\dot x = SW_{0}(x)Y(x)P^{-1}\hat{Z}(x)$ is
\begin{align}\label{dotV:thm}
\dot{V}(x)
%&= -\hat{Z}(x)^{\top} P^{-1} \big( \eta_{1}(x)+\eta_{2}(x) \big) P^{-1}\hat{Z}(x)\\
& = -\hat{Z}(x)^{\top} P^{-1}
\begin{bmatrix}
I_{p} \\ S^{\top}\frac{\partial\hat{Z}(x)}{\partial x}^{\top}
\end{bmatrix}^{\top}\notag\\
&\quad\cdot\begin{bmatrix}
0_{p\times p} & -Y(x)^{\top}W_{0}(x)^{\top}\\
- W_{0}(x) Y(x) & 0_{(q+N)\times p}
\end{bmatrix}\notag\\
&\quad\cdot\begin{bmatrix}
I_{p} \\ S^{\top}\frac{\partial\hat{Z}(x)}{\partial x}^{\top}
\end{bmatrix}
P^{-1}\hat{Z}(x).
\end{align}
On the other hand, since we choose a $\hat{Z}(x)$ whose first $n$ components coincide with $x$, $\frac{\partial\hat{Z}(x)}{\partial x}$ has full column rank. Hence, the quadratic constraint in Assumption \ref{assumption:boundaryD} can be equivalently written as
\begin{align*}
\frac{\partial\hat{Z}(x)}{\partial x}D_{0}D_{0}^{\top}\frac{\partial\hat{Z}(x)}{\partial x}^{\top} - \frac{\partial\hat{Z}(x)}{\partial x}R_{D}R_{D}^{\top}\frac{\partial\hat{Z}(x)}{\partial x}^{\top} \preceq 0
\end{align*}
and then expressed as
\begin{equation} \label{condition:assumptionD0}
\begin{array}{l}
\begin{bmatrix}
I_{p} \\ D_{0}^{\top}\frac{\partial\hat{Z}(x)}{\partial x}^{\top}
\end{bmatrix}^\top \\
\cdot
\begin{bmatrix}
- \frac{\partial\hat{Z}(x)}{\partial x} R_{D} R_{D} ^{\top}
\frac{\partial\hat{Z}(x)}{\partial x}^{\top} & 0_{N\times T}\\
0_{T\times N} & I_T
\end{bmatrix}
\begin{bmatrix}
I_{p} \\ D_{0}^{\top}\frac{\partial\hat{Z}(x)}{\partial x}^{\top}
\end{bmatrix}
\preceq 0.
\end{array}
\end{equation}
Using the equation $D_{0} = X_{1} - S \overline{W}_{0}$ which follows
from \eqref{data-equation} with
\begin{align}
\overline{W}_{0}:=
\begin{bmatrix}
\,\overline{U}_{0} \\ \,Z_{0}
\end{bmatrix}\in \mathbb{R}^{(q+N)\times T},
\end{align}
we can derive the following condition equivalent to \eqref{condition:assumptionD0}
\begin{align}\label{condition:assumptionD}
&\begin{bmatrix}
I_{p} \\ S^{\top}\frac{\partial\hat{Z}(x)}{\partial x}^{\top}
\end{bmatrix}^{\top}\notag\\
&\quad\cdot\begin{bmatrix}
\frac{\partial\hat{Z}(x)}{\partial x}\big( X_{1}X_{1}^{\top}-R_{D}R_{D}^{\top} \big) \frac{\partial\hat{Z}(x)}{\partial x}^{\top} & -\frac{\partial\hat{Z}(x)}{\partial x}X_{1}\overline{W}_{0}^{\top} \\
-\overline{W}_{0}X_{1}^{\top} \frac{\partial\hat{Z}(x)}{\partial x}^{\top} & \overline{W}_{0}\overline{W}_{0}^{\top}
\end{bmatrix}\notag\\
&\quad\cdot\begin{bmatrix}
I_{p} \\ S^{\top}\frac{\partial\hat{Z}(x)}{\partial x}^{\top}
\end{bmatrix}\preceq 0.
\end{align}

%\cdpcomment{Note to bear in mind during the numerical simulations: condition  \eqref{condition:assumptionD} suggests to choose $R_D$ such that
%$X_{1}X_{1}^{\top}-R_{D}R_{D}^{\top}\succeq 0$. This is consistent with the TAC paper, where we were setting $R_D = X_1$.}

%\cdpcomment{Remark for a better understanding: the matrix in \eqref{dotV:thm} can be obtained from the matrix \eqref{matrix:boundaryE} by replacing $E$ with $S$ using the identity
%$X_{1} -E \hat{U}_{0}(x)= SW_{0}\cdp{(x)}$. The same replacement in the matrix in
%\eqref{matrix:boundaryE} returns the matrix in \eqref{condition:assumptionD}.
%NO need to include this remark. }
%
%{\color{olive} Replacing $E$ with $S$ in
%\eqref{matrix:boundaryE} won't return the matrix in \eqref{condition:assumptionD}. There will be an additional matrix $\text{diag}(I_{q},0_{N})$ in the lower right corner.}
%\medskip

Using these conditions, we present a data-driven stabilizing design free of Assumption \ref{assumption:boundaryB} as follows:

\begin{theorem}\label{noisy-stabilization-II}
For the polynomial system (\ref{polysyst:original}), under Assumption \ref{assumption:boundaryD}, if there exist a
%symmetric
positive definite matrix $P\in\mathbb{R}^{p\times p}$, matrix polynomial $Y(x)\in\mathcal{P}^{T\times p}$, $\epsilon_{1}(x)\in\Sigma$, and $\epsilon_{2}(x)\in\Sigma$ such that
\begin{align}
&Z_{0}Y(x) = H(x)P,\label{condition:thm-1}\\
%&P-\mu I_{p} \succeq 0,\label{condition:thm-2}\\
& \begin{bmatrix}
\Upsilon_{D}(x)\!- \epsilon_{1}(x)I_{p}  & -Y(x)\!^{\top}\!W_{0}(x)\!^{\top}\!-\epsilon_{2}(x) \frac{\partial\hat{Z}(x)}{\partial x}X_{1}\overline{W}_{0}^{\top} \\
%- W_{0}(x) Y(x) -\epsilon_{2}(x) \overline{W}_{0}X_{1}^{\top} \Frac{\partial\hat{Z}(x)}{\partial x}^{\top} &
\star &
\epsilon_{2}(x)\overline{W}_{0}\overline{W}_{0}^{\top}
\end{bmatrix}\notag\\
&\quad\quad\quad\quad\quad\quad\quad\quad\quad\quad\quad\quad\quad\quad\quad\quad\quad \in\Sigma^{p+q+N},\label{condition:thm-3}
\end{align}
where %$\epsilon_{1}(x)\in\Sigma$, $\epsilon_{2}(x)\in\Sigma$ for all $x\in\mathbb{R}^{n}$, and
\begin{align}\label{UpsD}
\Upsilon_{D}(x) = \epsilon_{2}(x) \frac{\partial\hat{Z}(x)}{\partial x}\big( X_{1}X_{1}^{\top}-R_{D}R_{D}^{\top} \big) \frac{\partial\hat{Z}(x)}{\partial x}^{\top},
\end{align}
then the state-feedback controller \eqref{control:thm1}
%\begin{align}\label{control:thm}
%u = U_{0}Y(x)P^{-1}\hat{Z}(x)
%\end{align}
%globally
stabilizes the polynomial system. Moreover, if $\epsilon_{1}(x)>0$ for all $x\ne 0$, the closed-loop system is globally asymptotically stable.
\end{theorem}
%
%\cdpcomment{It seems to me that conditions \eqref{condition:thm-1}-\eqref{condition:thm-3} are equivalent to conditions \eqref{condition:prop-1}-\eqref{condition:prop-3}. In that respect, Proposition 1 and Theorem 1 are equivalent. Do you agree? On the other hand, I don't see the equivalence between Theorem 1 and Corollary 1. See my comment below. }
%
%\pt{As for the comment by Claudio, I don't see the
%equivalence because in Theorem 1 we don't have the
%assumption on $B$.}
%
%{\color{olive}Proposition 1 and Theorem 1 are not equivalent in my opinion. First, Assumption \ref{assumption:boundaryB} is not required in Theorem 1. Second, replacing $E$ with $S$ does not change \eqref{matrix:boundaryE} into \eqref{condition:assumptionD}.}

\proof
Under condition (\ref{condition:thm-3}), it holds that
\begin{align*}
&\begin{bmatrix}
\Upsilon_{D}(x) & -Y(x)^{\top}W_{0}(x)^{\top}\!-\epsilon_{2}(x) \frac{\partial\hat{Z}(x)}{\partial x}X_{1}\overline{W}_{0}^{\top} \\
%- W_{0}(x) Y(x) -\epsilon_{2}(x) \overline{W}_{0}X_{1}^{\top} \Frac{\partial\hat{Z}(x)}{\partial x}^{\top} &
\star &
\epsilon_{2}(x)\overline{W}_{0}\overline{W}_{0}^{\top}
\end{bmatrix}\\ \succeq &
\begin{bmatrix}
\epsilon_{1}(x)I_{p}  & 0_{p\times(q+T)} \\
0_{(q+T)\times p} & 0_{(q+T)\times(q+T)}
\end{bmatrix}.
\end{align*}
Recall the expression of $\dot{V}$ (\ref{dotV:thm})
along the closed-loop system and the expression (\ref{condition:assumptionD})
of Assumption \ref{assumption:boundaryD}.
Then, for any $\epsilon_{2}(x)\ge0$ $\forall x\in\mathbb{R}^{n}$, it holds that
\begin{align}%\label{dotV:thm}
\dot{V}(x) &\le
-\hat{Z}(x)^{\top} P^{-1}
\begin{bmatrix}
I_{p} \\ S^{\top}\frac{\partial\hat{Z}(x)}{\partial x}^{\top}
\end{bmatrix}^{\top}\notag\\
&\quad \cdot
\begin{bmatrix}
%\epsilon_{2}(x) \frac{\partial\hat{Z}(x)}{\partial x}\big( X_{1}X_{1}^{\top}\!-\!R_{D}R_{D}^{\top} \big) \frac{\partial\hat{Z}(x)}{\partial x}^{\top}
\Upsilon_{D}(x)
& -Y(x)^{\top}W_{0}(x)^{\top}\!-\!\epsilon_{2}(x) \frac{\partial\hat{Z}(x)}{\partial x}X_{1}\overline{W}_{0}^{\top} \\
%- W_{0}(x) Y(x)\! -\!\epsilon_{2}(x) \overline{W}_{0}X_{1}^{\top} \frac{\partial\hat{Z}(x)}{\partial x}^{\top}
\star
&\epsilon_{2}(x)\overline{W}_{0}\overline{W}_{0}^{\top}
\end{bmatrix}\notag\\
&\quad \cdot \begin{bmatrix}
I_{p} \\ S^{\top}\frac{\partial\hat{Z}(x)}{\partial x}^{\top}
\end{bmatrix}
P^{-1}\hat{Z}(x).
\end{align}
Using this inequality and condition (\ref{condition:thm-3}), we can conclude that
\begin{align*}
\dot{V}(x) \le  -\epsilon_{1}(x)\hat{Z}(x)^{\top} P^{-1}\cdot P^{-1}\hat{Z}(x)~~\forall x\in\mathbb{R}^{n},
\end{align*}
and the thesis follows.
\qed
\medskip

Theorem \ref{noisy-stabilization-I} and Theorem \ref{noisy-stabilization-II} use the same stabilizing approach but different data parameterizations of the closed-loop system. Comparing these two results, Theorem \ref{noisy-stabilization-II} is more advantageous in a few aspects. First, by pulling out the unknown dynamics $A$ and $B$ in one matrix $S$, Theorem \ref{noisy-stabilization-II} no longer requires Assumption \ref{assumption:boundaryB}. Second, the SOS condition \eqref{condition:thm-3} can be made more computationally efficient by a change of decision variable, as we discuss in the next subsection. Nonetheless, Theorem \ref{noisy-stabilization-I} is important because it provides insights on how to arrive at Theorem \ref{noisy-stabilization-II} while establishing connections with the robust stabilization results of \cite{cdpTAC2020}.

\begin{remark}{\rm (A least-square-based design)}
%In Subsection \ref{subsection:noisy2}  this identity will be used to bring about a variation of Theorem \ref{noisy-stabilization-I} with some distinguishing features.
%Here we use it to establish a link between Theorems \ref{noisy-stabilization-I} and \ref{noisy-stabilization-II} and a sequential system-identification-plus-robust-control kind of result.
%In the results above we have dealt with the uncertainty on the model induced by the noisy data by exploiting Assumption \ref{assumption:boundaryD} on the matrix of noise. An alternative to this approach is to start from the identity \eqref{data-equation} and adopt an estimate of the model that minimizes the Frobenius norm of $X_1-S\overline W_0=D_0$, which returns the least square estimate $S_*:=X_1 \overline W_0^\dag$ of $S$.
Identity \eqref{data-equation} can be used to establish a variation of Theorem \ref{noisy-stabilization-I} that has a similar feature as Theorem \ref{noisy-stabilization-II} but is based on a system's representation that uses a least-square estimate of the system's matrices.
By \eqref{data-equation}, we can adopt an estimate of the model $S$ that minimizes the Frobenius norm of $X_1-S\overline W_0=D_0$. This choice returns the least-square estimate $S_*:=X_1 \overline W_0^\dag$ of $S$, with $\overline W_0^\dag$ the pseudo-inverse of $\overline W_0$,  and entails an estimation error $S_*-S$.
We note that the estimate $S_*$ can be used for a model-based control design in which the actual model $S$ is replaced by the estimate  $S_*$, with the expectation that for smaller $S_*-S$, the result will be increasingly accurate. To show this, let us consider %for the sake of simplicity
the case of a full-row rank matrix $\overline W_0$, so that $S_*-S$ is given by $D_0 \overline W_0^\dag$. Then, we obtain the %additional
representation
\begin{align*}
X_{1} -E \hat{U}_{0}(x)= SW_{0}(x)= S_* W_{0}(x) - D_0 \overline W_0^\dag W_{0}(x),
\end{align*}
Based on this representation and on the same arguments of the proof of Theorem \ref{noisy-stabilization-I}, one realizes that, under Assumption \ref{assumption:boundaryD} and without  Assumption \ref{assumption:boundaryB}, an analogous of Theorem \ref{noisy-stabilization-I} holds provided that in condition \eqref{condition:I-3} the matrices $X_1, \hat{U}_{0}(x)$ are replaced by $S_* W_{0}(x), \overline W_0^\dag W_{0}(x)$ respectively.
\end{remark}

\subsection{Computationally more efficient stabilization conditions}
%Stabilization conditions independent of sampling number $T$
\label{subsection:noisy3}

Recalling that $F(x)= U_{0}Y(x)P^{-1}$, we can obtain $F(x)P= U_{0}Y(x)$. Defining $K(x):=F(x)P$, we can write the term $W_{0}(x)Y(x)$ in \eqref{dotV:thm} as
\begin{align*}
&\quad~ W_{0}(x)Y(x) \\
&=
\begin{bmatrix}
W(x)U_{0}Y(x) \\ Z_{0}Y(x)
\end{bmatrix}=
\begin{bmatrix}
W(x)F(x)P \\ H(x)P
\end{bmatrix}=
\begin{bmatrix}
W(x)K(x) \\ H(x)P
\end{bmatrix},
\end{align*}
which suggests to adopt $K(x)$ as a new decision variable, modify condition \eqref{condition:thm-3} and obtain the following result:

\begin{corollary}\label{noisy-stabilization-III}
For the polynomial system (\ref{polysyst:original}), under Assumption \ref{assumption:boundaryD}, if there exist a positive definite matrix
%{\color{blue}
$P\in\mathbb{R}^{p\times p}$, matrix polynomial $K(x)\in\mathcal{P}^{m\times p}$, $\epsilon_{1}(x)\in\Sigma$, and $\epsilon_{2}(x)\in\Sigma$
%}
such that
\begin{align}
%&P-\mu I_{p} \succeq 0,\label{condition:crlry-1}\\
& \begin{bmatrix}
\Upsilon_{D}\!(x)\! -\!\epsilon_{1}(x)I_{p}
& %-[K(x)^{\top}W(x)^{\top}~~PH(x)^{\top}]
-\!\begin{bmatrix}
W\!(x)K(x) \\ H(x)P
\end{bmatrix}^{\top}\!
-\!\epsilon_{2}(x)\frac{\partial\hat{Z}(x)}{\partial x}X_{1}\!\overline{W}_{0}^{\top}\\
\star & \epsilon_{2}(x) \overline{W}_{0} \overline{W}_{0}^{\top}
\end{bmatrix}\notag\\
&\quad\quad\quad\quad\quad\quad\quad\quad\quad\quad\quad\quad\quad\quad\quad\quad\quad
\in\Sigma^{p+q+N},\label{condition:crlry-2}
\end{align}
where %$\epsilon_{1}(x)\in\Sigma$,
%$\epsilon_{2}(x)>0$ for all $x\in\mathbb{R}^{n}$,
%$\epsilon_{2}(x)\in\Sigma$,
%and
$\Upsilon_{D}(x)$ is as in \eqref{UpsD},
%{\color{violet}
%\begin{align*}
%\Upsilon_{D}(x) = \epsilon_{2}(x) \frac{\partial\hat{Z}(x)}{\partial x} \big(X_{1}X_{1}^{\top}- R_{D} R_{D} ^{\top} \big)\frac{\partial\hat{Z}(x)}{\partial x}^{\top},
%\end{align*}
%}
then the state-feedback controller
\begin{align}\label{control:crlry}
u = K(x)P^{-1}\hat{Z}(x)
\end{align}
%globally
stabilizes the polynomial system. Moreover, if $\epsilon_{1}(x)>0$ for all $x\ne 0$, the closed-loop system is globally asymptotically stable.
\end{corollary}

\proof
Under the control law (\ref{control:crlry}), the closed-loop system can be expressed as
\begin{align*}
\dot{x}  = S
\begin{bmatrix}
W(x)K(x)P^{-1} \\ H(x)
\end{bmatrix}\hat{Z}(x).
\end{align*}
Then, the time derivative of the Lyapunov function $V(x)=\hat{Z}(x)^{\top}P^{-1}\hat{Z}(x)$ is
\begin{align*}
 \dot{V}(x) &= \hat{Z}(x)^{\top}P^{-1} \Bigg( \frac{\partial\hat{Z}(x)}{\partial x}S
 \begin{bmatrix}
W(x)K(x) \\ H(x)P
\end{bmatrix}\\
&\quad\quad +
\begin{bmatrix}
W(x)K(x) \\ H(x)P
\end{bmatrix}^{\top} S^{\top}\frac{\partial\hat{Z}(x)}{\partial x}^{\top} \Bigg)
P^{-1}\hat{Z}(x)\\
& = -\hat{Z}(x)^{\top} P^{-1}
\begin{bmatrix}
I_{p} \\ S^{\top}\frac{\partial\hat{Z}(x)}{\partial x}^{\top}
\end{bmatrix}^{\top}\\
&\quad\cdot
\begin{bmatrix}
0_{p\times p} & -K(x)^{\top}W(x)^{\top} & -PH(x)^{\top}\\
-W(x)K(x) & 0_{q\times q} & 0_{q\times N} \\
-H(x)P  & 0_{N\times q} & 0_{N\times N}
\end{bmatrix}\\
&\quad\cdot
\begin{bmatrix}
I_{p} \\ S^{\top}\frac{\partial\hat{Z}(x)}{\partial x}^{\top}
\end{bmatrix}
P^{-1}\hat{Z}(x).
\end{align*}
Using the expression (\ref{condition:assumptionD}) of Assumption \ref{assumption:boundaryD} and the condition (\ref{condition:crlry-2}), we can conclude that
\begin{align*}
\dot{V}(x) \le -\epsilon_{1}(x)\hat{Z}(x)^{\top} P^{-1}P^{-1}\hat{Z}(x)~~\forall x\in\mathbb{R}^{n}.
\end{align*}
The thesis follows.
%Hence, the closed-loop system is stable, and it is globally stable if $\|\hat{Z}(x)\|$ is radially unbounded. Moreover, if $\epsilon_{1}(x)>0$ for all $x\ne 0$, then $\dot{V}(x)<0$ for all $x\ne 0$, and thus the closed-loop system is globally asymptotically stable at the origin.
\qed

\medskip

In Theorem \ref{noisy-stabilization-II}, the size of decision variable $Y(x)$ is dependent on the number of sampled data $T$. For systems having higher order and degree, $T$ can be large which will result in heavy computational burden in the SOS programming.
%\cdpcomment{Why $T$ must be large? From Remark \ref{rem:rank_Z_0} it doesn't seem so. Is it because we are missing accurate bounds on $T$?}
Hence, by changing the decision variable $Y(x)$ used in Theorem \ref{noisy-stabilization-II}, the SOS condition \eqref{condition:thm-3} becomes \eqref{condition:crlry-2}, which is independent of $T$, and the stabilizer design is more computationally efficient.

%\cdpcomment{I added the remark showing that the condition
%\eqref{condition:crlry-2} is an SOS relaxation of a pointwise matrix S-lemma condition, which is a pointwise necessary and sufficient condition. In the red text I make comparisons with Szanier's paper. Let's continue to think what is the best form of this comparison that is beneficial to the paper.}

\begin{remark}
{\rm (An SOS relaxation of a pointwise necessary and sufficient condition)}
%We elaborate here on an interesting interpretation  of the previous results.
%%Theorem \ref{thm.noisy.stab} and Corollary \ref{cor.noisy.stab.II}.
Observe that investigating the stabilization problem via the Lyapunov function $V(x)=\hat{Z}(x)^\top P^{-1} \hat{Z}(x)$ and the feedback controller $u=K(x)P^{-1}\hat{Z}(x)$, the problem becomes one of finding the matrix $P\succ 0$ and the  matrix $K(x)$ such that for all $x\in \mathbb{R}^n\setminus \{0\}$, the matrix in
%\eqref{Vdot.general}
the proof of Corollary \ref{noisy-stabilization-III}
satisfies
\[
\begin{bmatrix}
0_{p\times p} & -K(x)^{\top}W(x)^{\top} & -PH(x)^{\top}\\
-W(x)K(x) & 0_{q\times q} & 0_{q\times N} \\
-H(x)P  & 0_{N\times q} & 0_{N\times N}
\end{bmatrix}
\succ 0
\]
for all $S$ such that \eqref{condition:assumptionD} holds.
%\[
%\begin{bmatrix}
%I_N \\
%S^\top \frac{\partial Z}{\partial x}^\top
%\end{bmatrix}
%^\top
%\begin{bmatrix}
%R(x) R(x)^\top - \frac{\partial Z}{\partial x} X_1 X_1^\top \frac{\partial Z}{\partial x}^\top & \frac{\partial Z}{\partial x} X_1 W_0^\top\\[2mm]
%W_0 X_1^\top \frac{\partial Z}{\partial x}^\top  & -W_0 W_0^\top
%\end{bmatrix}
%\begin{bmatrix}
%I_N \\
%S^\top \frac{\partial Z}{\partial x}^\top
%\end{bmatrix}
%\succeq 0.
%\]
For each fixed $x\in \mathbb{R}^n\setminus \{0\}$, one can pointwise apply the matrix S-lemma of \cite{vanwaarde2020noisy} to this formulation. Under the technical condition that for each  $x\in \mathbb{R}^n\setminus \{0\}$, there exists a matrix $\overline S$ for which \eqref{condition:assumptionD} with $S$ replaced by $\overline S$ holds,  the stabilization problem is solvable if and only if there exist a matrix $P\succ 0$ and a  matrix $K(x)$ such that  for all $x\in \mathbb{R}^n\setminus \{0\}$   there exist $\epsilon_1(x)> 0$ and $\epsilon_2(x)\ge 0$ for which
the matrix \eqref{condition:crlry-2} is positive semi-definite. Having this condition satisfied
%for all $x\in \mathbb{R}^n\setminus \{0\}$
leads to an infinite dimensional problem. A natural way to overcome this obstacle is to relax the positive semi-definiteness condition by requiring $K(x), \epsilon_1(x), \epsilon_2(x)$ to be matrix polynomials
%polynomial matrices
and the resulting matrix polynomial \eqref{condition:crlry-2} to be an SOS matrix, which is the condition in Corollary \ref{noisy-stabilization-III}.
%
%\cdpcomment{words in red above to be refined: $\epsilon_2(x)$ is our statement must be strictly positive for all $x$, while it can be $=0$ as a consequence of the $S$-lemma.}
%
%{\color{violet} Theorem 1 needs $\epsilon_{2}(x)>0$ as the matrix-valued Young's inequality is used. In Theorem 2 and Corollary 1, instead of Young's inequality, I think we indeed use the S-lemma in a pointwise manner. Thus, the constraint on $\epsilon_{2}(x)$ should be $\epsilon_{2}(x)\ge 0$ in Theorem 2 and Corollary 1.}
\end{remark}

\begin{remark}
{\rm (Another convex relaxation)}
The stabilization problem of polynomial systems from data has been tackled in \cite{dai2020semialgebraic} using density functions and polyhedral constraints on the uncertainties, which leads to Quadratically Constrained Quadratic Program whose convex relaxation is solved via moments-based techniques. Stemming from the results in \cite{cdpTAC2020,guoCDC2020}, our approach uses Lyapunov second stability theorem and quadratic constraints, leading to %substantially different
SOS programs whose solutions directly provide stabilizing controllers and Lyapunov functions.
%The two results  differ substantially because we use direct Lyapunov stability theory and quadratic constraints, while  \cite{dai2020semialgebraic} adopts %Rantzer's
%an
%approach to stability based on density functions and polyhedral constraints on the uncertainties  %in our approach
%As a result, the  stability condition \eqref{condition:crlry-2} is an easily checkable SOS program only dependent on the size of the system's matrices while
%\cite{dai2020semialgebraic} obtains a Quadratically Constrained Quadratic Program dependent on the size of the data, whose convex relaxation is solved via moments-based techniques.
\end{remark}

%\cdpcomment{In the write-up I have a remark (Remark 5) showing that the condition
%\eqref{condition:crlry-2} is an SOS relaxation of a pointwise matrix S-lemma condition, which is a pointwise necessary and sufficient condition. Szanier's paper on polynomial systems also uses a pointwise necessary and sufficient condition, but his approach is based on Rantzer's dual approach, linear (and not quadratic) stability conditions and Farkas Lemma, which makes the analysis much less clear. We should discuss whether or not the addition of this remark,  as well as the addition of a comparison with Sznaier's results, is beneficial to the paper.}

%\cdpcomment{In general some discussion should be added to the mere derivation of results and make the paper ``juicer". We should lure the reader to appreciate the paper.}

%\cdpcomment{Numerical simulations should be added possibly comparing Proposition \ref{noisy-stabilization-I} and Corollary \ref{noisy-stabilization-III}, although it is intuitive to expect a better performance of the latter.}

%\cdpcomment{We should discuss whether, beside stabilization with noisy data, additional results should be investigated. For instance, we could revisit the bilinear result and give a local stabilization result with guaranteed basin of attraction for polynomial systems. Is it worthwhile? The optimal control problem for polynomial systems is likely related to the noisy LQR problem and I am not sure we are ready to deal with the noisy optimal polynomial stabilization case.}

%\bigskip

\section{Example}
\label{section:examples}
In this section, we apply the data-driven control method to the stabilization of the Van der Pol oscillator, which is a popular benchmark to test data-driven control results. The simulations are conducted using the SOSTOOLS in MATLAB.

%\subsection{Van der Pol oscillator}
Consider the controlled Van der Pol oscillator
\begin{align}\label{VanderPol}
\dot{x}_{1} &= x_{2},\notag\\
\dot{x}_{2} &= -x_{1}+(1-x_{1}^2)x_{2}+u.
\end{align}

Let $Z(x)$ be the power vector that contains all monomials of $x$ having degrees $1$ to $3$. Choose $\hat{Z}(x)=[x_{1}~~x_{2}]^{\top}$, then
\begin{align*}
H(x)=\begin{bmatrix}
1 & 0 & x_2 & x_1 & 0 & x_1x_2 & 0 & x_1^2 & 0 \\
0 & 1 & 0 & 0 & x_2 & 0 & x_1x_2 & 0 & x_2
\end{bmatrix}^{\top}.
\end{align*}
Run an experiment with initial condition $x_{0}=[-0.1~~0.1]^{\top}$ and input $u=\sin t$ from $t=0$ to $t=10$. The sampling period is $0.5$.

The data $X_{1}$ is corrupted by a noise
$D_{0}= 0.05 X_{1}$. We set $T=12$ and $R_{D}R_{D}^{\top}=0.1X_{1}X_{1}^{\top}$. SOS programs are then formulated to search for $P$, $\epsilon_{2}$, and the coefficients of $\epsilon_{1}(x)$ and $Y(x)$ with both $\epsilon_{1}(x)$ and $Y(x)$ having degree $2$.

First, we apply the special case of Theorem \ref{noisy-stabilization-I} where the input vector field is independent of $x$ (Remark \ref{rem:inp-state-ind}). The SOSTOOLS obtains the solution
\begin{align*}
  P&=\begin{bmatrix}
   8.67\times 10^{-7} & -2.199\times 10^{-7}\\
 -2.199\times 10^{-7} &   1.925\times 10^{-7}
  \end{bmatrix},\\
    \epsilon_{1}(x)&=x_{1}(5.203\times 10^{-5}x_{1} + 1.627\times 10^{-6}x_{2})\\
   &\quad+ x_{2}(1.627\times 10^{-6}x_{1} + 5.39\times 10^{-5}x_{2}),\\
 \epsilon_{2}&=5.71\times 10^{-7}.
\end{align*}
The stabilizer is calculated as
\begin{align*}
u&= - x_{1} (9.2 x_{1}^2 +2.2 x_{1} x_{2} + x_{1} + 8.1 x_{2}^2  +0.24 x_{2} + 9.7)\\
&\quad - x_{2} (38 x_{1}^2 +6.1 x_{1} x_{2} + 5.5 x_{1} + 33 x_{2}^2 + 2.3 x_{2} + 29).
\end{align*}

Using the design conditions in Theorem \ref{noisy-stabilization-II}, the solution given by SOSTOOLS is
\begin{align*}
  P&=\begin{bmatrix}
  1.692\times 10^{-6} &  1.311\times 10^{-7} \\
  1.311\times 10^{-7} & 1.104\times 10^{-6}
   \end{bmatrix},\\
     \epsilon_{1}(x)&=x_{1}(9.165\times 10^{-6}x_{1} + 6.389\times 10^{-10}x_{2})\\
   &\quad+ x_{2}(6.389\times 10^{-10}x_{1} + 9.258\times 10^{-6}x_{2}),\\
  \epsilon_{2}&=1.061\times 10^{-6}.
\end{align*}
The designed stabilizer is
\begin{align*}
  u&= -  x_{1} (3.7\times 10^{-4} x_{1}^2 +5.1\times 10^{-5} x_{1} x_{2} - 0.021 x_{1}\\
   &\quad\quad\quad +2.2\times 10^{-4} x_{2}^2 -0.092 x_{2} + 0.94)\\
 &\quad -  x_{2} (-1\times 10^{-4} x_{1}^2 -5.3\times 10^{-3} x_{1} x_{2} - 3.1\times 10^{-3} x_{1} \\
  &\quad\quad\quad +5.5\times 10^{-4} x_{2}^2 -0.064 x_{2} + 1.2).
  %u&= -  x_{1} (7.5\times 10^{-3} x_{1}^2 -9.1\times 10^{-3} x_{1} x_{2} - 0.036 x_{1}\\
%   &\quad\quad\quad -1.4\times 10^{-3} x_{2}^2 +6.3\times 10^{-3} x_{2} + 1.3)\\
% &\quad -  x_{2} (9.4\times 10^{-3} x_{1}^2 +2.2\times 10^{-3} x_{1} x_{2} - 1.2\times 10^{-3} x_{1} \\
%  &\quad\quad\quad -0.02 x_{2}^2 +0.012 x_{2} + 3.6).
\end{align*}
%\cdpcomment{Perhaps it is useful to add also the Lyapunov function. Sznaier can't compute the Lyapunov function -- I think -- and having a Lyapunov function is useful for other studies.}
%The closed-loop system is asymptotically stable at the origin.

Alternatively, we can set $K(x)$ having degree $2$ as a decision variable and utilize Corollary \ref{noisy-stabilization-III} to design the data-driven stabilizer. The solution to the SOS program is
\begin{align*}
  P&=\begin{bmatrix}
   6.332\times 10^{-7} &  -5.787\times 10^{-8}\\
  -5.787\times 10^{-8} & 6.934\times 10^{-7}
%   4.486\times 10^{-6} &  3.048\times 10^{-7}\\
%  3.048\times 10^{-7} & 1.87\times 10^{-6}
  \end{bmatrix},\\
  \epsilon_{1}(x)&=x_{1}(5.125\times 10^{-6}x_{1} + 1.366\times 10^{-10}x_{2})\\
     &\quad+ x_{2}(1.366\times 10^{-10}x_{1} + 5.16\times 10^{-6}x_{2}),\\
%  \epsilon_{1}(x)&=x_{1}(2.704\times 10^{-5}x_{1} + 2.034\times 10^{-9}x_{2})\\
%   &\quad+ x_{2}(2.034\times 10^{-9}x_{1} + 2.718\times 10^{-5}x_{2}),\\
  \epsilon_{2}&=7.184\times 10^{-7}.
%  \epsilon_{2}(x)&=2.492\times 10^{-6}.
\end{align*}
The stabilizer is
\begin{align*}
  u&= -x_{1} (1.9\times 10^{-4} x_{1}^2  -8.9\times 10^{-5} x_{1} x_{2} - 0.026 x_{1} \\
  &\quad\quad\quad +4.2 \times 10^{-4} x_{2}^2 -0.091 x_{2} +1.5) \\
  &\quad - x_{2} (4.8\times 10^{-5} x_{1}^2  - 1.1\times 10^{-4} x_{1} x_{2} - 3.6\times 10^{-3} x_{1} \\
  &\quad\quad\quad - 2.4\times 10^{-4} x_{2}^2 -0.072 x_{2} +1.4).
%  u&= -x_{1} (1.3\times 10^{-4} x_{1}^2  -1\times 10^{-4} x_{1} x_{2} - 0.045 x_{1} \\
%  &\quad\quad\quad +2.3 \times 10^{-4} x_{2}^2 +0.01 x_{2} +1.2) \\
%  &\quad - x_{2} (1.3\times 10^{-4} x_{1}^2  - 2.9\times 10^{-4} x_{1} x_{2} - 8.8\times 10^{-3} x_{1} \\
%  &\quad\quad\quad + 4.3\times 10^{-4} x_{2}^2 +0.016 x_{2} +3.1).
\end{align*}

\begin{figure}
  \centering
  \includegraphics[width=0.4\textwidth]{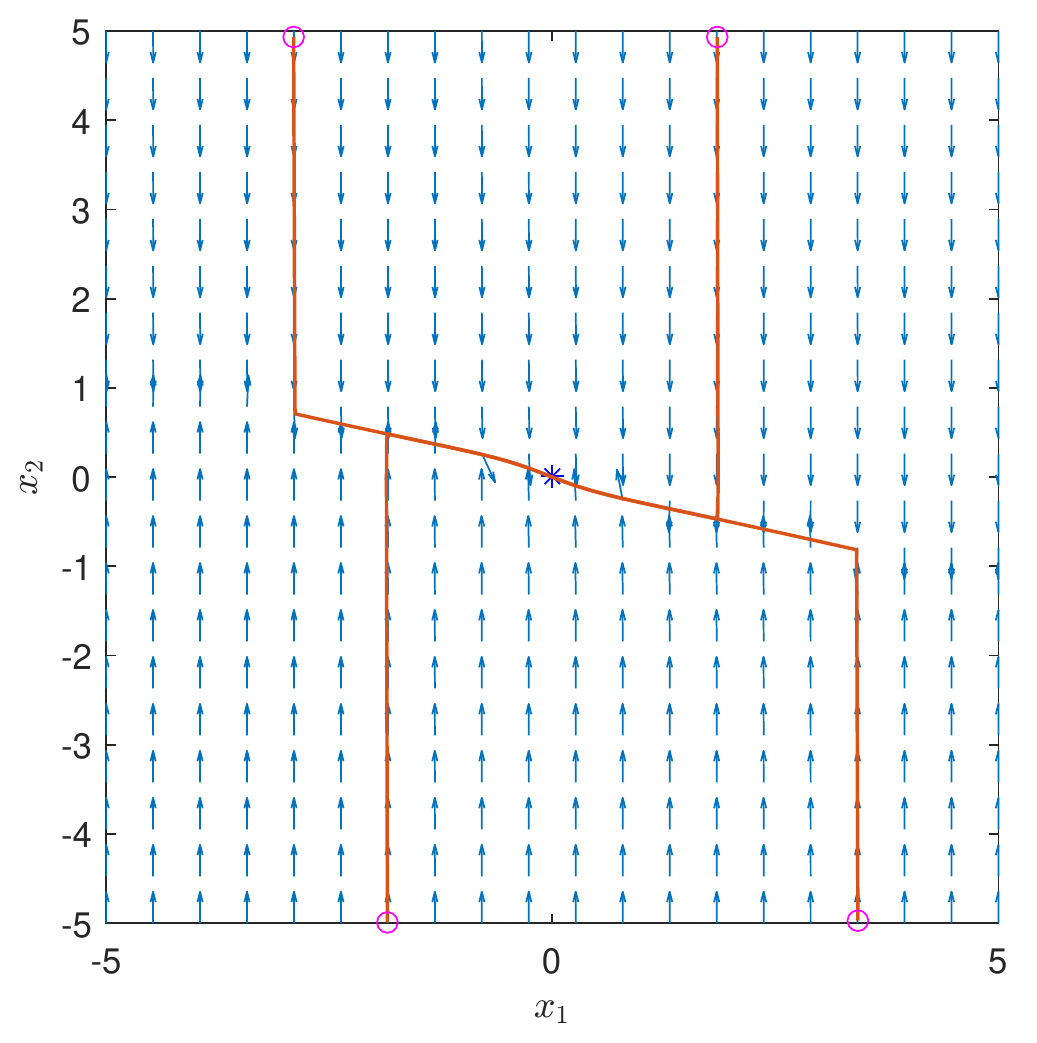}\\
  \caption{Phase portrait of the closed-loop system under the controller designed using Theorem \ref{noisy-stabilization-I}.}\label{figure:oscillatorPhPT1}
\end{figure}

\begin{figure}
  \centering
  \includegraphics[width=0.4\textwidth]{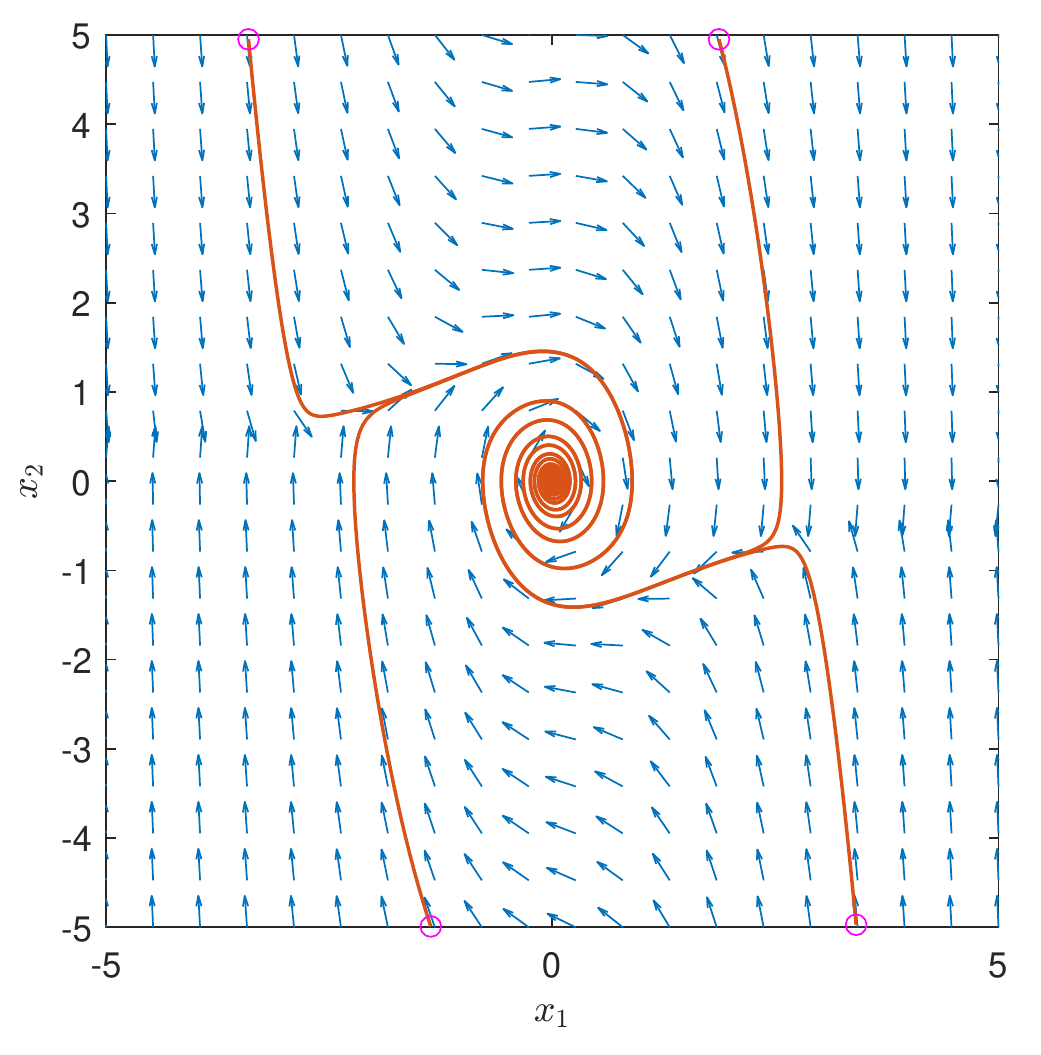}\\
  \caption{Phase portrait of the closed-loop system under the controller designed using Theorem \ref{noisy-stabilization-II}.}\label{figure:oscillatorPhPT2}
\end{figure}

\begin{figure}
  \centering
  \includegraphics[width=0.4\textwidth]{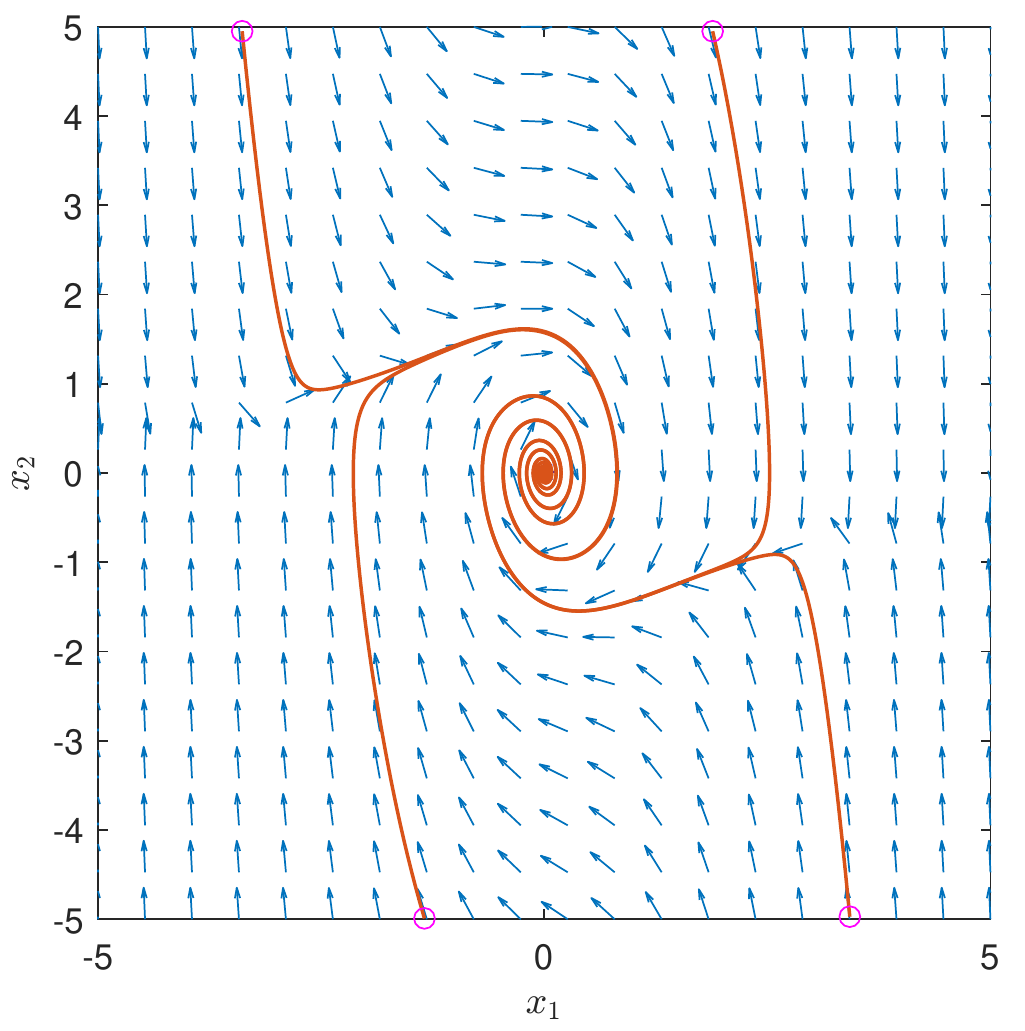}\\
  \caption{Phase portrait of the closed-loop system under the controller designed using Corollary \ref{noisy-stabilization-III}.}\label{figure:oscillatorPhPC1}
\end{figure}

The phase portraits of the closed-loop systems under the designed controllers are illustrated in Figures \ref{figure:oscillatorPhPT1} to \ref{figure:oscillatorPhPC1} respectively. The figures show that, designed using the same data set, the controllers stabilize system \eqref{VanderPol} at the origin with different transient performances. The computational time needed for formulating and solving the SOS program and then obtaining the control gain $F(x)$ is $146.5957s$ for Theorem \ref{noisy-stabilization-I}, $151.3056s$ for Theorem \ref{noisy-stabilization-II}, and $10.6662s$ for Corollary \ref{noisy-stabilization-III}. This verifies that Corollary \ref{noisy-stabilization-III} is more computational efficient than Theorems \ref{noisy-stabilization-I} and \ref{noisy-stabilization-II}.

%Again using Corollary \ref{noisy-stabilization-III}, if we set $T=7$ and $R_{D}R_{D}^{\top}=0.01X_{1}X_{1}^{\top}$, the solution is
%\begin{align*}
%  P&=\begin{bmatrix}
%   4.288\times 10^{-6} &  -1.339\times 10^{-6}\\
%  -1.339\times 10^{-6} & 5.285\times 10^{-6}
%  \end{bmatrix},\\
%  \epsilon_{1}(x)&=x_{1}(2.344\times 10^{-5}x_{1} -6.735\times 10^{-8}x_{2})\\
%   &\quad+ x_{2}(-6.735\times 10^{-8}x_{1} + 2.368\times 10^{-5}x_{2}),\\
%  \epsilon_{2}(x)&=8.149\times 10^{-5}.
%\end{align*}
%The stabilizer is
%\begin{align*}
%  u&= -x_{1} (6.2\times 10^{-3} x_{1}^2  +2.2\times 10^{-3} x_{1} x_{2} -0.047 x_{1} \\
%  &\quad\quad\quad +4.9 \times 10^{-3} x_{2}^2-0.016 x_{2} +2.1) \\
%  &\quad - x_{2} (-9.3\times 10^{-4} x_{1}^2  - 9.4\times 10^{-4} x_{1} x_{2} - 0.054 x_{1} \\
%  &\quad\quad\quad - 3.1\times 10^{-4} x_{2}^2 -0.094 x_{2} +1.5).
%\end{align*}
%}

\vspace{-0.45cm}
\section{Conclusion}
\label{section:conclusion}

We have shown how to synthesize stabilizers and Lyapunov functions for unknown nonlinear polynomial systems starting from noisy data and using Lyapunov second theorem. We do not assume the unknown noise to have any specific form as long as it has a known quadratic bound over the experiment. The state-dependent stabilizing gain is solved via SOS programs that are computationally tractable. Interestingly, the efficiency of the stabilization design can be improved by changing the data parameterization of the closed-loop system even when the same design method is utilized.

Nonlinear data-driven stabilization is a fundamental and important problem that lays the foundation for our ongoing works on nonlinear polynomial systems such as local control with guaranteed domain of attraction and optimal control with quadratic costs. Another interesting topic is looking into the computational aspects of SOS programming and further improving the efficiency of the data-driven designs. Our results can play an important role in learning control policies for those manifold applications where polynomial systems and SOS optimization have found wide use.

\def\bibfont{\small}
\bibliographystyle{IEEEtran}
\bibliography{referenceDD}
\end{document}